%% file: root.tex
\documentclass[letterpaper, 10 pt, conference]{ieeeconf}  % Comment this line out if you need a4paper

\IEEEoverridecommandlockouts                              % This command is only needed if 
\usepackage{graphics} % for pdf, bitmapped graphics files
\usepackage{amsmath} % assumes amsmath package installed
\usepackage{amssymb}  % assumes amsmath package installed
\usepackage{graphicx}
\usepackage{algorithm}
\usepackage{algorithmic}
\usepackage{comment}

\usepackage[margin=1in]{geometry}

\usepackage[section]{placeins}

\def\R{\mathbb{R}}

\usepackage{xcolor}
\usepackage{multirow}

\allowdisplaybreaks

\title{\LARGE \bf
Iteratively Preconditioned Gradient-Descent Approach \\
for Moving Horizon Estimation Problems
}

\author{Tianchen Liu$^{1}$, Kushal Chakrabarti$^{2}$, Nikhil Chopra$^{1}$ 
\thanks{$^{1}$ Department of Mechanical Engineering, University of Maryland, College Park, MD 20742, USA. {\tt\small \{tianchen, nchopra\}@umd.edu}} %
\thanks{$^{2}$ Division of Data \& Decision Sciences, Tata Consultancy Services Research, Mumbai, 400607 India. {\tt\small chakrabarti.k@tcs.com}} %
}

\begin{document}

\maketitle
\thispagestyle{empty}
\pagestyle{empty}

\input{sections/00_abstract.tex}

\input{sections/01_intro.tex}

\input{sections/03_problem.tex}

\input{sections/04_approach}

\input{sections/05_experiments.tex}

\input{sections/06_conclusion.tex}

\input{sections/appendix.tex}

\bibliographystyle{IEEEtran}
\bibliography{refs_mhe}

\end{document}

%% file: sections/00_abstract.tex
\begin{abstract}
Moving horizon estimation (MHE) is a widely studied state estimation approach in several practical applications. In the MHE problem, the state estimates are obtained via the solution of an approximated nonlinear optimization problem. However, this optimization step is known to be computationally complex. Given this limitation, this paper investigates the idea of iteratively preconditioned gradient-descent (IPG) to solve MHE problem with the aim of an improved performance than the existing solution techniques. To our knowledge, the preconditioning technique is used for the first time in this paper to reduce the computational cost and accelerate the crucial optimization step for MHE. The convergence guarantee of the proposed iterative approach for a class of MHE problems is presented. Additionally, sufficient conditions for the MHE problem to be convex are also derived. Finally, the proposed method is implemented on a unicycle localization example. The simulation results demonstrate that the proposed approach can achieve better accuracy with reduced computational costs.
\end{abstract}

%% file: sections/01_intro.tex
\section{Introduction}
Moving horizon estimation (MHE), which first appeared in~\cite{michalska1995moving}, is an optimization-based solution approach for state estimation of linear and nonlinear systems. This approach models and solves an optimization problem at each sampling instant for obtaining the best state variable values in a finite period, regarded as state estimates. While using complete prior information for the estimation should generate better estimates; however, the computation cost can quickly become intractable as well. MHE handles this challenge by utilizing a finite number of past measurements and control inputs and discarding the previous information to maintain a feasible computational cost when solving the optimization problem. Compared with other state estimation methods, e.g., Kalman filter~\cite{kalman1960new}, extended Kalman filter (EKF)~\cite{gelb1974applied} and particle filter~\cite{gordon1993novel}, MHE performs well for the constrained state estimation problem when the arrival cost is accurately approximated, which contains information of the discarded data. A general introduction and some applications of MHE can be found in~\cite{alessandri2020moving}. Due to its performance and efficiency, MHE has become a widely used approach for state estimation for many applications, including linear (e.g.,~\cite{rao2001constrained,alessandri2003receding,alessandri2008moving}) and nonlinear (e.g.,~\cite{alessandri2005robust,alessandri2010advances,osman2021generic,schiller2022suboptimal}) systems. Stability analysis has also been investigated for specific scenarios (e.g.,~\cite{liu2013moving,ji2015robust,muller2017nonlinear}).

As discussed in~\cite{zou2020moving}, the performance of MHE critically relies on the algorithm used to solve the underlying optimization problem. Various strategies have been developed to reduce computational complexity while maintaining estimation accuracy. In~\cite{morabito2015simple}, Nesterov’s fast gradient method expedites the optimization step. However, the limitation is that the approach only applies to linear systems. In~\cite{hashemian2015fast}, nonlinear system equations are approximated by Carleman linearization expressions to reduce the computational cost of obtaining the gradient and Hessian. In~\cite{alessandri2017fast}, three approaches based on the gradient, conjugate gradient, and Newton's method have been proposed to reduce the computational effort using an iterative approach, which are demonstrated to be more effective than the Kalman filter by simulation. The proposed approach in this paper is similar to the one using Newton's method in~\cite{alessandri2017fast}. Still, instead of calculating the inverse of the Hessian matrix, we introduce a preconditioner matrix~\cite{kelley1999iterative} into the iterative approach. This technique can significantly help to avoid the costly step of inverting a matrix.

% \kushal{The transition from here to our work in the next paragraph is abrupt. Maybe conclude the above paragraph or add a new one with why we need a new approach for solving MHE, following the literature review above, and what needs to be improved compared to them. Add a sentence about it in the abstract also.} \tianchen{I have rephrased the paragraph and the abstract}

Notably, we transform the distributed iteratively preconditioned gradient-descent (IPG) approach proposed in~\cite{chakrabarti2021accelerating} to its centralized counterpart and employed it for the nonlinear state estimation in the MHE framework. To our knowledge, this is the first attempt to utilize the preconditioning technique to reduce computational costs and accelerate the optimization step for MHE. The convergence proof of the proposed algorithm and the sufficient conditions for the MHE problem to be convex are presented. To our knowledge, such convexity analysis has not been found in previous literature. Finally, the proposed approach is implemented on a simulation example for estimating the locations of a unicycle mobile robot. The results demonstrate that the MHE approach achieves better performance than EKF, invariant EKF (InEKF)~\cite{barrau2016invariant}, and a recently developed IPG observer~\cite{chak2023obsv}. Compared with the default optimization solver in Matlab, the proposed approach can obtain the same results with a reduced computational cost. The main contributions of this paper are summarized below.
\begin{itemize}
    \item An iterative optimization method using the preconditioning technique, following~\cite{chakrabarti2021accelerating}, is employed to solve the MHE problem. A convergence proof of the proposed approach is presented.
    \item Since the above convergence guarantee requires convexity, sufficient conditions for the formulated MHE problem to be convex are derived.
    \item The proposed approach is implemented on a mobile robot localization problem, and the results are compared with other nonlinear state estimators to demonstrate the advantages.
\end{itemize}

The rest of the paper is organized as follows. The MHE problem is described in Section II. Section III proposes a solution to the MHE problem utilizing the IPG method. Empirical results are presented in Section IV to illustrate the performance of the proposed solution compared to other state estimation methods. Finally, Section V concludes this paper.

%% file: sections/03_problem.tex
\section{Problem Description}
In this paper, we consider the discrete-time system
\begin{align}
	% x_{i+1} & = f(x_i, u_i) + \varepsilon_{x_i} \label{eq:sys_dynamics} \\
	% y_i & = h(x_i) + \varepsilon_{y_i} \label{eq:sys_obs} \qquad (i = 0, \dots, T)
    x_{i+1} & = f(x_i, u_i) + w_i, \label{eq:sys_dynamics} \\
	y_i & = h(x_i) + v_i, \label{eq:sys_obs} \qquad (i = 0, \dots, T-1) %\\
    %w_i & \in \mathbb{W}, \quad v_i \in \mathbb{V} \nonumber
\end{align}
where $x_i \in \mathbb{R}^n$, $u_i \in \mathbb{R}^m$, and $y_i \in \mathbb{R}^p$ denote the states, the inputs, and the observations at $i^{th}$ sampling instant, respectively. The process disturbance set $\mathbb{W} \subseteq \mathbb{R}^n$ and measurement noise set $\mathbb{V} \subseteq \mathbb{R}^p$ are assumed to be compact
with $0 \in \mathbb{W}$ and $0 \in \mathbb{V}$~\cite{rao2003constrained}. Hence, $w_i \in \mathbb{W}$ and $v_i \in \mathbb{V}$ are bounded process disturbance and measurement noise vectors. The system drift function $f: (\mathbb{R}^n,\mathbb{R}^m) \rightarrow \mathbb{R}^n$ and the measurement function $h: \mathbb{R}^n \rightarrow \mathbb{R}^p$ are assumed to be known linear or nonlinear functions. $T$ is the total number of sampling steps. The formulation indicates that an input control and a (partial) measurement occur at each sampling step. 
% \kushal{Is Gaussian required?} \tianchen{Gaussian not required, but boundedness is required. Boundedness requirement does not influence the proof of proposition 2. I will change the example part accordingly.}

%The idea of the moving horizon estimation is established on the nonlinear optimization of the disturbances of the system. 
% In this paper, we consider the case when the process disturbances and measurement noises are additive. 
For the sampling instant $t = N, \dots, T$, the MHE problem can be formulated as the following nonlinear optimization problem~\cite{rawlings2006particle},
\begin{align}
\min_{x_{\{t-N:t\}}} & \Phi_{t-N} := \sum_{i=t-N}^{t-1} \Big( w_i^T {Q^{-1}} w_i + v_i^T {R^{-1}} v_i \Big)  \nonumber \\
& \qquad \qquad + \Gamma(x_{t-N}), \nonumber \\
\text{s.t. } w_i & = x_{i+1} - f(x_{i}, u_i), \label{eq:mhe_formulation} \\
 v_i & = y_i - h(x_{i}). \quad (i=t-N,\ldots,t-1) \nonumber
\end{align}
% \kushal{if the above formulation exists in prior works, we should cite one of them} \tianchen{Thanks, I have added the reference here.}
The optimization variables are $x_{t-N}$ to $x_t$, given the past $N$ known control inputs and observations. $Q \in \mathbb{R}^{n \times n}$ and $R \in \mathbb{R}^{p \times p}$ are diagonal positive definite weighting matrices for the disturbances. The arrival cost $\Gamma(x_{t-N})$ summarizes the discarded past information. We employ an EKF-based approximation of the arrival cost~\cite{rao2003constrained},
\begin{align}
    \Gamma(x_{t-N}) = & \, (x_{t-N} - \hat{x}_{t-N})^T \Pi_{(t-N)}^{-1} (x_{t-N} - \hat{x}_{t-N}) \nonumber \\
    & + \Phi_{t-N}^*.
\end{align}
where $\hat{x}_{t-N}$ 
% \kushal{This is different than $\hat{x}_{t-N}$ in~(10), right? If true, then maybe we change the notation from hat to something else here. } \tianchen{This and $\hat{x}_{t-N}$ in~(10) have the same meaning, i.e., the estimate value of ${x}_{t-N}$ from the optimization step. The subtle difference is that it's actually $\hat{x}_{t-N+1}$ obtained in~(10) that is used for solving $\Gamma(x_{t-N+1})$ and $\Phi_{t+1}$ for the next time step. Please feel free to revise the text to clarify this point.} 
and $\Phi_{t-N}^*$ are the estimate of ${x}_{t-N}$ and the optimal objective function value obtained at the previous time instant. 
% \kushal{by optimizing $\phi_{t-N}$ or $\phi_{t-1}$?} 
% \tianchen{By optimizing $\Phi_{t-N}$. I revised the notation $\Phi_{t-N}$ in~(3). It's clearer this way.} 
$\Pi_{(t-N)} \in \mathbb{R}^{n \times n}$ is a positive definite weighting matrix, updated for the next time instant via the following matrix Riccati equation~\cite{rao2003constrained},
\begin{align}
S_1 & = J_h \Pi_{(t-N)} J_h^T + R, \nonumber \\
S_2 & = J_f \Pi_{(t-N)} J_h^T (S_1)^{-1} J_h \Pi_{(t-N)} J_f^T, \nonumber \\
\Pi_{(t-N+1)} & = J_f \Pi_{(t-N)} J_f^T - S_2 + Q, 
\label{eq:update_Pi_i}
\end{align}
where $J_f \in \mathbb{R}^{n \times n}$ and $J_h \in \mathbb{R}^{p \times n}$ are the Jacobian of $f$ and $h$ with respect to state variables $x$. 

To present our results, we require a few more notations.
For each $t \geq N$, $Y^{(t)} \in \mathbb{R}^{Np}$ and $U^{(t)} \in \mathbb{R}^{Nm}$ denote the concatenating column vectors of the past $N$ consecutive measurements and control inputs before $t^{th}$ time instant, respectively, i.e.,  $Y^{(t)} = \left[y_{t-N}^T, \dots, y_{t-1}^T \right]^T$ and $U^{(t)} = \left[u_{t-N}^T, \dots, u_{t-1}^T \right]^T$.
%\begin{align}
%    U^{(t)} = & \left[u_{t-N}^T, \dots, u_{t-1}^T \right]^T \\
%    Y^{(t)} = & \left[y_{t-N}^T, \dots, y_{t-1}^T \right]^T \label{eqn:meas_vec}
%\end{align}
We let $\| \cdot \|$, $\lambda_{max} [\cdot]$, and $\lambda_{min} [\cdot]$ denote the induced 2-norm, the largest, and smallest eigenvalue of a matrix. 

%For each $i \geq N$, an estimate of $x_i$ is obtained by solving for $x_{i-N+1}$ in~\eqref{eqn:meas_vec} and propagating the obtained solution of $x_{i-N+1}$ forward by $N$ sampling instants using the dynamics function $f$. In this paper, we solve this observer problem using the IPG approach~\cite{chak2023obsv}.

%For the problem as shown in Eq.~\ref{eq:mhe_formulation}, we can verify that the $F$ is twice continuously differentiable if the arrival cost $\Gamma(x_{t-N})$ is quadratic with respect to $x_{t-N}$, and minimum of $F$ exists and is finite.

%% file: sections/04_approach.tex
\section{Proposed Approach}
In this section, we present details of the proposed IPG approach, henceforth referred to as MHE-IPG, employed to solve the MHE problem. The convergence proof of MHE-IPG is provided. Finally, a convexity analysis of the MHE problem is discussed. 

\subsection{MHE-IPG Approach}
The main contribution of the proposed IPG approach lies in utilizing a preconditioning technique to accelerate the optimization step when solving MHE problem. For each time instant $t=N,\ldots, T$, the steps of MHE-IPG are as follows.

\textbf{Step 1}. For $t^{th}$ time instant, define the optimization variable vector $\xi^{(t)} \in \mathbb{R}^{(N+1)n}$, which is the concatenating column vector of variables $x_{t-N}, \dots, x_{t}$, i.e.,
\begin{equation}
    \xi^{(t)} = [x_{t-N}^T, \dots, x_{t}^T]^T.
\end{equation}
Then, the MHE problem at time instant $t \geq N$ in Eq. (\ref{eq:mhe_formulation}) is equivalent to
\begin{align}
& \min_{\xi^{(t)}} F(\xi^{(t)}, U^{(t)}, Y^{(t)}) \nonumber \\
= & \sum_{i=t-N}^{t-1} (x_{i+1} - f(x_{i}, u_i))^T {Q^{-1}} (x_{i+1} - f(x_{i}, u_i)) \nonumber \\
& + \sum_{i=t-N}^{t-1} (y_i - h(x_{i}))^T {R^{-1}} (y_i - h(x_{i}))  \nonumber \\
& + (x_{t-N} - \hat{x}_{t-N})^T \Pi_{(t-N)}^{-1} (x_{t-N} - \hat{x}_{t-N}).
\label{eq:converted_mhe}
\end{align}
For $t=N$, the initial state estimate $\hat{x}_{0}$ and a positive definite matrix $\Pi_{(0)}$ are chosen. Otherwise, $\hat{x}_{t-N}$ and $\Pi_{(t-N)}$ are obtained from estimates at the previous time instant.

\textbf{Step 2}. This step solves the optimization problem~\eqref{eq:converted_mhe} using the idea of IPG. IPG is an iterative algorithm. At each iteration $k = 0,1,\ldots$, an estimate $\xi_k^{(t)}$ and a preconditioner matrix $K \in \mathbb{R}^{(N+1)n \times (N+1)n}$ are maintained. Before the iterations start, we select the positive scalar constants $\beta, \delta$, and initialize $\xi^{(t)}_0, K_{0}$. At each iteration $k$, the current estimate and preconditioner are updated to
\begin{align}
\xi^{(t)}_{k+1} & = \xi_k^{(t)} - \delta K_{k} g(\xi^{(t)}_k),\label{eq:update_x}\\
K_{k+1} & = K_{k} - \alpha_k [(H(\xi_k^{(t)}) + \beta I) K_{k} - I],\label{eq:update_K}
\end{align}
until $\| \xi^{(t)}_{k+1} - \xi^{(t)}_k \|< \epsilon$, where $\epsilon$ is a predetermined tolerance value. $g(\xi_k^{(t)})$ and $H(\xi_k^{(t)})$ denote the gradient and Hessian of $F$ with respect to $\xi$, evaluated at $\xi = \xi_k^{(t)}$. $I$ is the identity matrix with the same dimension as $H(\xi_k^{(t)})$ and $K_{k}$. $\alpha_k$ is selected following the condition in Eq.~(\ref{eq:select_alpha_ipg}). Let $\hat{\xi}^{(t)} = \xi^{(t)}_{k+1}$ and go to Step 3.

\textbf{Step 3}. Decompose $\hat{\xi}^{(t)}$ to the state estimates for the past $N$ time instants as
\begin{equation}
    \hat{\xi}^{(t)} = [\hat{x}_{t-N}^T, \dots, \hat{x}_{t}^T]^T.
\end{equation}
Record $\hat{x}_{t-N}^T, \dots, \hat{x}_{t}^T$ to the MHE results of (\ref{eq:converted_mhe}) and go to the next time instant $t+1$ to estimate $\xi^{(t+1)}$. Utilizing $\hat{\xi}^{(t)}$, $\Pi_{(t-N+1)}$ is updated via Eq.~(\ref{eq:update_Pi_i}). Moreover, we form the initial estimate $\xi^{(t+1)}_0$ for next time instant as
\begin{equation}
    \xi^{(t+1)}_0 = [\hat{x}_{t-N+1}^T, \dots, \hat{x}_{t}^T, f(\hat{x}_{t}^T)]^T.
\end{equation}
Go back to Step 1 for solving (\ref{eq:converted_mhe}) for the time instant $t+1$ using $\xi^{(t+1)}_0$, $\hat{x}_{t-N+1}$ and $\Pi_{(t-N+1)}$. Repeat until the estimation is done for all sampling steps. 

Note that the MHE-IPG approach can be extended to solve constrained MHE problems. Due to space constraints, this will be considered in future work.

\subsection{Convergence Analysis of MHE-IPG}
We make the following assumptions to present our convergence results of the proposed approach. \\
\textbf{Assumption 1.} The system equations $f$ and $h$ are assumed to satisfy certain conditions such that $F(\xi)$ is convex and twice continuously differentiable, with the minimum solution(s) of~\eqref{eq:converted_mhe} exist and denoted as $\xi^{(t)*} \in \Xi^{(t)*}$. For brevity, we will denote $\xi^{(t)*}$ as $\xi^{*}$.\\
\textbf{Assumption 2.} The Hessian of $F(\xi)$, denoted by $H(\xi)$, is assumed to be Lipschitz continuous with respect to the 2-norm with Lipschitz constant $\gamma$, i.e.,
$$\|H(\xi_1) - H(\xi_2)\| \leq \gamma \| \xi_1-\xi_2 \|, \, \forall \xi_1, \xi_2 \in \R^{(N+1)n}.$$
We further assume that $\| H(\xi^*) \|$ is upper bounded as $ \| H(\xi^*) \| \leq q$ for some $q \in (0,\infty)$ and $H(\xi^*)$ is non-singular at any minimum point $\xi^{*} \in \Xi^{*}$. \\
\textbf{Assumption 3.} The gradient of $F(\xi)$, denoted by $g(\xi)$, is assumed to be Lipschitz continuous with a positive Lipschitz constant $l$, i.e.,
$$\|g(\xi_1) - g(\xi_2)\| \leq l \| \xi_1-\xi_2 \|, \, \forall \xi_1, \xi_2 \in \R^{(N+1)n}.$$

For each time instant $t$, we introduce the following notation. We define the `optimal' preconditioner matrix
\begin{equation}
    K^{*} = (H(\xi^*)+\beta I)^{-1}.
\end{equation}
It can be concluded that $K^{*}$ is well-defined with $\beta > 0$ and Assumption 1. We let $\eta$ denote its induced $2$-norm
\begin{align}
    \eta = & \|K^{*}\| = \|(H(\xi^*) + \beta I)^{-1}\| = \frac{1}{\lambda_{min}[H(\xi^*)] + \beta}. \nonumber 
\end{align}
For each iteration $k$, we define
\begin{align}
    \Tilde{K}_k & = K_k - K^*,
    \label{eq:error_K_k}
\end{align}
the coefficient for convergence of $K_k$,
\begin{equation}
    \rho_k = \|I - \alpha_k (H(\xi_k)+\beta I) \|,
\end{equation}
and the estimation error $z_k = \xi_k - \xi^*$. Let $\rho = \sup \rho_k$. If $\beta > 0$ and $0 < \alpha_k < \frac{1}{\lambda_{max} [H(\xi_k)] + \beta}$,
%\begin{equation}
%    0 < \alpha_k < \frac{1}{\lambda_{max} [H(\xi_k)] + \beta}
%    \label{eq:select_alpha}
%\end{equation}
then we conclude that $\rho_k \in [0, 1), \forall k \geq 0$ (see~\cite{chakrabarti2021accelerating}, Lemma 1). 

The following lemma is essential for the convergence of our proposed method.

\textbf{Lemma 1.~\cite{chakrabarti2021accelerating}}
For each time instant $t \geq N$, consider the IPG update~\eqref{eq:update_x}-\eqref{eq:update_K} with parameters $\beta, \delta, \alpha_k > 0$. Then, under Assumptions~1-3, $\|\Tilde{K}_{k}\|$ is bounded by
\begin{align}
    \|\Tilde{K}_{k+1}\| \leq & \rho^{k+1} \|\Tilde{K}_0\| + \gamma \eta (\alpha_k \|z_k\|  \nonumber \\
    & + \rho \alpha_{k-1} \|z_{k-1}\| + \cdots + \rho^k \alpha_{0} \|z_{0}\|)
    \label{eq:K_upper_bound_2}
\end{align} 
%$\lambda_{min}[H(\xi^*)] > 0$, $\eta < \frac{1}{\beta}$.
The detailed proof can be found in~\cite{chakrabarti2021accelerating} Appendix A.3.

Next, we present the convergence result of the proposed IPG approach for solving the MHE problem~\eqref{eq:converted_mhe} for any $t \geq N$. We drop the superscript $(t)$ from the presented results for brevity.

\textbf{Proposition 1.} 
Suppose that Assumptions~1-3 holds. For each time instant $t \geq N$, consider the IPG update~\eqref{eq:update_x}-\eqref{eq:update_K} with parameters $\beta > 0$, $\delta > 0$, %$\delta \in (0,1]$
and $ \alpha_k \in (0, \frac{1}{\lambda_{max} [H(\xi_k)] + \beta})$. Let the initial estimate $\xi_0$ and preconditioner matrix $K_{0}$ be selected to satisfy 
\begin{equation}
    %\frac{\eta \gamma}{2} \|x_0 - \xi^*\| + l \|K_0 - K^{*}\| + \eta \beta \leq \frac{1}{2\mu}
    \frac{\delta \eta \gamma}{2} \|\xi_0 - \xi^*\| + \eta \beta + \eta q | 1-\delta | + \delta l \|K_{0} - K^{*}\| \leq \frac{1}{2\mu},
    \label{eq:init_select_cond}
\end{equation}
where $\mu \in (1, \frac{1}{\rho})$ and $\eta = \|K^{*}\|$. If 
\begin{equation}
    \hspace{-0.2em} \alpha_k < \min \{\frac{1}{\lambda_{max} [H(\xi_k)] + \beta}, \frac{\mu^{k}(1-\mu \rho)}{2 l (1-(\mu \rho)^{k+1})} \}, \hspace{-0.2em}
    \label{eq:select_alpha_ipg}
\end{equation}
then for $k \geq 0$, 
\begin{equation}
    \|\xi_{k+1} - \xi^*\| < \frac{1}{\mu} \|\xi_k - \xi^*\|.
\end{equation}

The proof is deferred to Appendix A.

Proposition~1 implies that for each $t \geq N$, the estimates of the IPG approach locally converge to a solution of~\eqref{eq:converted_mhe} with a linear convergence rate of at least $\frac{1}{\mu}$.

%Since $\mu > 1$, Proposition 1 implies that the sequence $\xi_k$ locally converges to a solution $\xi^*$ with a linear convergence rate $\frac{1}{\mu}$. 

\begin{comment}
To obtain a simpler condition on $\alpha_k$,
compared to \ref{eq:select_alpha_ipg}, we can use a more conservative upper bound. Let $\Lambda$ be the upper bound of $\lambda_{max} [H(\xi_k)]$ for all $k$, and given $0 < \mu \rho < 1$, $\mu > 0$,
\begin{equation}
    \frac{\mu^k (1-\mu \rho)}{2l(1-(\mu \rho)^{k+1})} > \frac{\mu^k (1-\mu \rho)}{2l}
\end{equation}
we have
\begin{equation}
    \alpha_k < \min \{\frac{1}{\Lambda + \beta}, \frac{\mu^k (1-\mu \rho)}{2l}\}
    \label{eq:select_alpha_ipg_simplified}
\end{equation}
\end{comment}

%This is only true when the system dynamics $f(\cdot)$ and the state observer functions $h(\cdot)$ satisfy certain conditions. 

\subsection{MHE Convexity Analysis}
\label{sec:convexity}
Assumption 1 requires the converted function $F(\cdot)$ in~\eqref{eq:converted_mhe} to be convex and twice continuously differentiable. Hence, we present sufficient conditions on the system dynamics $f(\cdot)$ and observation function $h(\cdot)$ given known $U$ and $Y$ such that Assumption 1 holds.

Given the MHE formulation (\ref{eq:converted_mhe}), it is easy to observe that the Hessian $H(\xi) \in \mathbb{R}^{(N+1)n \times (N+1)n}$ with respect to $\xi$ is a tridiagonal block matrix, i.e., $H(\xi)$ is of the form
\begin{align}
H(\xi) & = {\footnotesize	\begin{bmatrix}
* & * &  &  &  &  \\
* & * & * &  &  &  \\
    & * & * & * &  & \\
    &   & \ddots & \ddots & \ddots &  \\
    &  &  & * & * & * \\
    &  &  &  & * & *
\end{bmatrix} }. \nonumber
\end{align}
where non-zeros elements exist only in partitioned blocks represented by $(*)$. Then we define the matrix $\overline{\overline{H}}_{i} \in \mathbb{R}^{2n \times 2n}$ as \\
\begin{align}
\overline{\overline{H}}_{i} & = \begin{bmatrix}
A_{11} & - J_{f}^T |_{x_i} Q^{-1} \\
- Q^{-1} J_{f}|_{x_i} & Q^{-1}
\end{bmatrix},
\label{eq:bar_bar_H_i} \\
\text{where } & A_{11} = \Pi^{-1} + \Tilde{J}_{f(0)} + V_{f(0)}^T \Tilde{Q} \Tilde{H}_{f(0)} + \Tilde{J}_{h(0)} \nonumber \\ 
& \quad \quad \quad + V_{h(0)}^T \Tilde{R} \Tilde{H}_{h(0)}, \quad (i = 0) \nonumber \\
&  A_{11} = \Tilde{J}_{f(i)} + V_{f(i)}^T \Tilde{Q} \Tilde{H}_{f(i)} + \Tilde{J}_{h(i)}   \nonumber \\ 
& \quad \quad \quad + V_{h(i)}^T \Tilde{R} \Tilde{H}_{h(i)}, \quad (i = 1, \dots, N-1) \nonumber
\end{align}
$(\cdot) |_{x_i}$ means to evaluate the expressions at the time instant $i$. The matrices are calculated as follows,
\begin{align}
    \Tilde{J}_{f(i)} & = J_f^T Q^{-1} J_f |_{x_i}, \quad \Tilde{J}_{h(i)} = J_h^T R^{-1} J_h |_{x_i}, \\
    V_{f(i)} & = I_{n} \otimes (f(x_i, u_i) - x_{i+1}), \\
    V_{h(i)} & = I_{n} \otimes (h(x_i) - y_{i}), \\
    \Tilde{Q} & = I_{n} \otimes Q^{-1}, \quad \Tilde{R} = I_{n} \otimes R^{-1}, \\
    \Tilde{H}_{f(i)} & = 
    \begin{bmatrix}
    \mathcal{H}_{f(1,1,1)} & \cdots & \mathcal{H}_{f(1,1,n)} \\
    \vdots &  & \vdots \\
    \mathcal{H}_{f(n,1,1)} & \cdots & \mathcal{H}_{f(n,1,n)} \\
    %\mathcal{H}_{f(1,2,1)} & \cdots & \mathcal{H}_{f(1,2,n)} \\
    %\vdots &  & \vdots \\
    %\mathcal{H}_{f(n,2,1)} & \cdots & \mathcal{H}_{f(n,2,n)} \\
    \vdots &  & \vdots \\
    \mathcal{H}_{f(1,n,1)} & \cdots & \mathcal{H}_{f(1,n,n)} \\
    \vdots &  & \vdots \\
    \mathcal{H}_{f(n,n,1)} & \cdots & \mathcal{H}_{f(n,n,n)} \\
    \end{bmatrix}_{|_{x_i}}, \\
    \Tilde{H}_{h(i)} & = \begin{bmatrix}
    \mathcal{H}_{h(1,1,1)} & \cdots & \mathcal{H}_{h(1,1,n)} \\
    \vdots &  & \vdots \\
    \mathcal{H}_{h(p,1,1)} & \cdots & \mathcal{H}_{h(p,1,n)} \\
    \vdots &  & \vdots \\
    \mathcal{H}_{h(1,n,1)} & \cdots & \mathcal{H}_{h(1,n,n)} \\
    \vdots &  & \vdots \\
    \mathcal{H}_{h(p,n,1)} & \cdots & \mathcal{H}_{h(p,n,n)} \\
    \end{bmatrix}_{|_{x_i}},
\end{align}
where $\otimes$ denotes the Kronecker product of two matrices. $\mathcal{H}_f \in \mathbb{R}^{n \times n \times n}$ and $\mathcal{H}_h \in \mathbb{R}^{p \times n \times n}$ are two 3-dimensional tensors concatenating Hessians of $f$ and $h$ with respect to $x$:
\begin{align}
    \mathcal{H}_{f(i,j,k)} = \frac{\partial^2 f_i}{\partial x_j \partial x_k}, \quad 
    \mathcal{H}_{h(i,j,k)} = \frac{\partial^2 h_i}{\partial x_j \partial x_k}. \nonumber 
\end{align}

\textbf{Proposition~2.} Consider the system dynamics function $f$ and the observation function $h$ in~\eqref{eq:sys_dynamics}-\eqref{eq:sys_obs}. If $\overline{\overline{H}}_{i}$, as defined in Eq.~(\ref{eq:bar_bar_H_i}), is positive semi-definite for all $i = 0, \dots, N-1$, then the MHE problem (\ref{eq:converted_mhe}) is convex.

The proof is deferred to Appendix B.

% \kushal{what is the dimension of $V_{f(i)}^T$ by definition~(22)?} \tianchen{The dimension of $V_{f(i)}$ is $\mathbb{R}^{n^2 \times n}$ and the dimension of $V_{h(i)}$ is $\mathbb{R}^{np \times n}$}

%\begin{comment}

Given Proposition 2, corollaries can be derived for the following special cases.

\textbf{Corollary 2.1.} Consider the system~\eqref{eq:sys_dynamics}-\eqref{eq:sys_obs} to be linear with state feedback control $u = -Kx$, i.e.,
\begin{equation}
    f(x) = Ax+Bu = (A-BK)x, \quad h(x) = Cx. \nonumber
\end{equation}
Then $J_f = A-BK =: A_c$, $J_h = C$, and $\Tilde{H}_{f(i)} = \Tilde{H}_{h(i)} = [0]$ for all $i$, and the sufficient condition~\eqref{eq:bar_bar_H_i} can be simplified as
\begin{equation}
\overline{\overline{H}} = \begin{bmatrix}
A_c^T Q^{-1} A_c + C^T R^{-1} C & - A_c^T Q^{-1} \\
- Q^{-1} A_c & Q^{-1}
\end{bmatrix} 
\end{equation}
being positive semi-definite. % for all $i = 0, \dots, N-1$.
%\end{comment}

\textbf{Definition 1.~\cite{horn2012matrix}} A square matrix $S$ is said to be diagonally dominant if
$|S_{ii}| \geq \sum_{j \neq i} |S_{ij}|, \, \forall i$,
i.e., for each row $i$ of the matrix, the absolute value of the diagonal element is no less than the sum of the absolute values of the rest elements in the same row. 

\textbf{Corollary 2.2.} 
Given $\overline{\overline{H}}_{i}$ and $A_{11}$ as defined in Eq.~(\ref{eq:bar_bar_H_i}), if $\overline{\overline{H}}_{i}$ is a diagonally dominant matrix and the diagonal entries of $A_{11}$ are non-negative for all $i = 0, \dots, N-1$, then the MHE problem (\ref{eq:converted_mhe}) is convex.
% \kushal{is $A_{11}$ symmetric?} \tianchen{Please check the added proof.}

\begin{proof}
First we show that $A_{11}$ is symmetric. It is easy to see that $\Pi^{-1}$, $\Tilde{J}_{f(i)}$ and $\Tilde{J}_{h(i)}$ are symmetric. Let the scalar $v_k$ denote the $k^{th}$ element of the vector $f(x_i, u_i) - x_{i+1}$, i.e., $f_k(x_i, u_i) - x_{i+1, k}$, and the scalar $q_k$ denote the entry $(k, k)$ of $Q^{-1}$. Then the entry $(i,j)$ of $V_{f(i)}^T \Tilde{Q} \Tilde{H}_{f(i)}$ is $\sum_{k=1}^{n} v_k q_k \mathcal{H}_{f(k,i,j)}$, and the entry $(j, i)$ is $\sum_{k=1}^{n} v_k q_k \mathcal{H}_{f(k,j,i)}$. Since $\mathcal{H}_{f(k,i,j)} = \mathcal{H}_{f(k,j,i)}$ for any $k$, $V_{f(i)}^T \Tilde{Q} \Tilde{H}_{f(i)}$ is symmetric. Similarly we can show that $V_{h(i)}^T \Tilde{R} \Tilde{H}_{h(i)}$ is symmetric. Hence, $A_{11}$ is symmetric, leading to a symmetric $\overline{\overline{H}}_{i}$.

A symmetric diagonally dominant matrix with real non-negative diagonal entries is positive semidefinite~\cite{horn2012matrix}. Hence, it can be concluded that as the above conditions are satisfied, then $H(\xi)$ is positive semidefinite, and the problem (\ref{eq:converted_mhe}) is convex.
\end{proof}

%% file: sections/05_experiments.tex
\section{Experiments}

In this section, the performance of the proposed MHE-IPG approach is evaluated on a location estimation problem of a mobile robot. We use the first-order Euler discretization to convert the continuous-time model formulation into a discrete-time one. The computations are performed in MATLAB 2022a on a Windows 11 laptop configured with an i7-9750H CPU and 16GB RAM.

The discrete-time kinematics of the unicycle model is formulated as follows,
\begin{align}
	x_{i+1} & = \begin{pmatrix}
		x_{i+1,1}\\
		x_{i+1,2} \\
		x_{i+1,3}
	\end{pmatrix} = \begin{pmatrix}
		f_1(x_{i}, u_{i}) + \varepsilon_{i, x_1}\\
		f_2(x_{i}, u_{i}) + \varepsilon_{i, x_2}\\
		f_3(x_{i}, u_{i}) + \varepsilon_{i, x_3}
	\end{pmatrix} \nonumber \\
    & = 
	\begin{pmatrix}
		x_{i,1} + dt \cdot u_{i, 1} \cos (x_{i,3}) + \varepsilon_{i, x_1}\\
		x_{i,2} + dt \cdot u_{i, 1} \sin (x_{i,3}) + \varepsilon_{i, x_2}\\
		x_{i,3} + dt \cdot u_{i, 2} + \varepsilon_{i, x_3}
	\end{pmatrix}.
\end{align}
State variables are the position $x_{i, 1}, x_{i, 2}$ and heading direction $x_{i, 3}$ in the world frame coordinates. Control inputs are $u_{i} = [u_{i, 1}, u_{i, 2}]^T$, where $u_{i, 1}$ is the forward speed and $u_{i, 2}$ is the angular velocity. $\varepsilon_{i,x}$ is the process disturbance vector.

The observations are the direct measurements of the position of the robot (e.g., global positioning system (GPS) measurements) with additive noises:
\begin{align}
	y_{i} = & \begin{pmatrix}
		y_{i, 1}\\
		y_{i, 2}
	\end{pmatrix}  = \begin{pmatrix}
		h_1(x_{i}, u_{i}) + \varepsilon_{i, y_1}\\
		h_2(x_{i}, u_{i}) + \varepsilon_{i, y_2}
	\end{pmatrix} = 
	\begin{pmatrix}
		x_{i, 1} + \varepsilon_{i, y_1} \\
		x_{i, 2} + \varepsilon_{i, y_2}
	\end{pmatrix}.
\end{align}
where $\varepsilon_{i, y}$ is the measurements noise vector.

The initial states are $x_0 = [0, 0, 0]^T$. The sampling time interval is $dt = 0.2$, and the total number of sampling instants is $T = 200$. The control inputs are
\begin{align}
	u_i & = 
	\begin{pmatrix}
		3 \\
		i/200
	\end{pmatrix}, \qquad \text{for $i = 0, ... , T-1.$}
\end{align}

The process noises $\varepsilon_{i, x_1}, \varepsilon_{i, x_2}, \varepsilon_{i, x_3} \sim N(0, 0.1)$, and the measurement noises $\varepsilon_{i, y_1}, \varepsilon_{i, y_2} \sim N(0, 0.4)$ are bounded with a maximal magnitude of $1.5$. {Given these parameters, it can be verified that Proposition 2 is valid for this problem.}

Different nonlinear estimators have been tested for this localization problem, including EKF, invariant EKF (InEKF)~\cite{barrau2016invariant}, IPG observer~\cite{chak2023obsv} and the MHE approach. EKF is a widely used technique for nonlinear state estimation but may suffer from divergence. InEKF, avoids the divergence issue by mapping the states to matrix Lie groups, where the converted problem is solved. The IPG observer was recently developed in~\cite{chak2023obsv} that uses the same iteratively preconditioning technique but in the manner of a Newton-type nonlinear observer. For the MHE approach, both the default Matlab unconstrained optimizer (`fminunc') and the proposed MHE-IPG algorithm are employed. 

To evaluate the performance of different state estimators, we use the root mean square error (RMSE) over $M$ simulation runs
\begin{equation}
    \bar{e} = \frac{1}{M} \sum_{m=1}^{M} \Big( \sum_{t=0}^{T} \|e_t^{(l)}\|^2 \Big)^{\frac{1}{2}}, %, t = 0, 1, \dots, T-1
\end{equation}
where $e_t^{(m)}$ is the estimation error of the $m$-th simulation run. $M=30$ runs are simulated with randomly generated noises.  

\begin{table}[ht!]
\centering
%\begin{tabular}{ |p{2cm}|p{2cm}|p{2cm}| }
\caption{Error Comparison for Different Estimators}
\label{table:error_comp}
\begin{tabular}{|c|c|c|c|}
%\hline
%\multicolumn{4}{|c|}{Time} \\
\hline
Method & Window Size $(N)$ & Mean $\bar{e}$ (m) & Variance $\bar{e}$ \\
\hline
Observations & - & 0.4966 &  0.0669 \\
EKF & 1 & 1.0626 &  0.6241 \\
InEKF & 1 &  0.1993 &  0.0150 \\
\hline
IPG Observer & 5 & 0.2809 &  0.0230 \\
MHE-default & 5 & 0.1943 &  0.0099 \\
MHE-IPG & 5 & 0.1943 &  0.0099 \\
\hline
IPG Observer & 10 & 0.2362 &  0.0178 \\
MHE-default & 10 & 0.1935 &  0.0104 \\
MHE-IPG & 10 & 0.1935 & 0.0104 \\
\hline
IPG Observer & 15 &  0.2462 &  0.0193 \\
MHE-default & 15 &  0.1867 &  0.0097 \\
MHE-IPG & 15 &  0.1867 &  0.0097 \\
\hline
IPG Observer & 20 & 0.4116 &  0.0308 \\
MHE-default & 20 &  0.1851 & 0.00956 \\
MHE-IPG & 20 & 0.1851 & 0.00956 \\
\hline
\end{tabular}
\end{table}

\begin{table}
\caption{Computational Cost for Different Estimators}
\label{table:compute_time}
\centering
\begin{tabular}{ |p{1.9cm}|p{1cm}|p{1.6cm}|p{1.6cm}|}
%\begin{tabular}{ cccc }
%\hline
%\multicolumn{4}{|c|}{Time} \\
\hline
Method & Window Size & Mean Computation Time (s) & Average Total Number of Iterations\\
\hline
EKF & 1 & 0.0009724 & - \\
InEKF & 1 & 0.010764 & -\\
\hline
IPG Observer & 5 & 9.9487 & -\\
MHE-default & 5 & 2.2322 & 3284.2 \\
MHE-IPG & 5 & 0.098335 & 2994.9\\
\hline
IPG Observer & 10 & 22.565 & -\\
MHE-default & 10 & 5.9298 & 4258.1 \\
MHE-IPG & 10 & 0.24038 & 3202.1\\
\hline
IPG Observer & 15 & 39.359 & - \\
MHE-default & 15 & 11.104 & 4486.3 \\
MHE-IPG & 15 & 0.57329 & 3034.4 \\
\hline
IPG Observer & 20 & 64.235 & - \\
MHE-default & 20 & 16.724 & 4459.1 \\
MHE-IPG & 20 & 0.82994 & 2863.1 \\
\hline
\end{tabular}
\end{table}

Table \ref{table:error_comp} shows the mean and variance of the estimation RMSE by different estimators against the ground truth. It can be seen that MHE results obtained by the `fminunc' solver and the proposed MHE-IPG approach outperform those of the other estimators regarding accuracy. In addition, the average RMSE of MHE results reduces as the window size $N$ increases, which makes sense since more information has been used for the estimation step. 

\begin{figure}[tbh!]
	\centering
	\includegraphics[width=\linewidth]{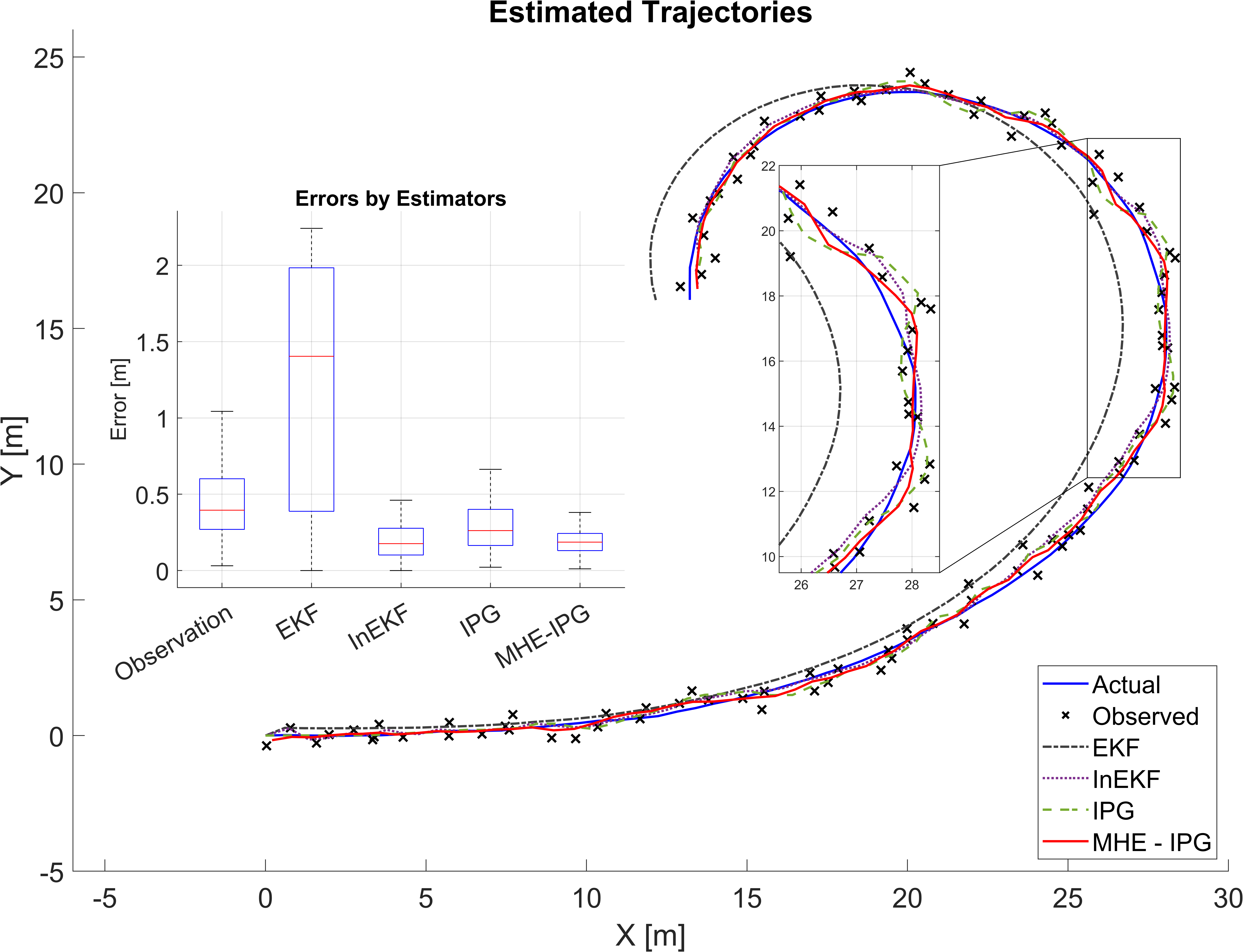}
	\caption[]{Estimated Trajectories and Box Plots for Errors ($N = 5$)}
	\label{fig:ex1_traj_n_5}
\end{figure}

Figure~\ref{fig:ex1_traj_n_5} shows the estimated trajectories and the box plots for estimation errors of one run when the window size $N=5$. The synthetic ground truth trajectory and the observations are plotted by the solid blue curve and black crosses, respectively. From the figure, it can be observed that EKF has relatively large errors. All other estimators can obtain lower mean and variance values of RMSE than the observations. Note that the IPG observer tends to be more influenced by the noisy observations, and thus has a slightly worse accuracy.

A main benefit of using the preconditioning technique is improving the optimization step's computation speed. Table \ref{table:compute_time} shows the average computational time for different estimators. The same stopping criteria is adopted to compare the proposed MHE-IPG algorithm against the default `fminunc' solver ($\epsilon = 10^{-6}$), and the average total numbers of iterations for the two solvers are listed. {$\beta = 0.5$ and $\delta = 1.6$ are used for the MHE-IPG solver.} It is observed that by employing the preconditioning technique, the proposed MHE-IPG algorithm can run faster and converge in fewer iterations for solving the MHE problem. For instance, the computational time for the window size $N=20$ by the MHE-IPG solver is lower than that by the default `fminunc' solver for the window size $N=5$. As mentioned above, a larger window size can lead to better estimation results. Hence, our proposed approach can be computationally less expensive to achieve the same level of accuracy or better with the same level of computational cost. {Furthermore, $\beta, \delta$ can be tuned to achieve faster convergence.}

%% file: sections/06_conclusion.tex
\section{Conclusion}
In this paper, we proposed a new iterative approach for solving the moving horizon estimation problem. To address the vital but computationally expensive optimization step in the MHE problem, the proposed MHE-IPG approach utilizes an iterative preconditioning technique. It helps in reducing the computational cost and accelerates the optimization step of MHE. Convergence of the proposed MHE-IPG approach is rigorously analyzed for convex cost functions and the sufficient conditions for the MHE problem to be convex are derived. Such conditions guarantee that the proposed MHE-IPG algorithm can obtain this solution, which has not been discussed in prior works, to the best of our knowledge. Finally, the simulated unicycle localization example highlights the improved accuracy of the MHE-IPG solution to other nonlinear state estimators, including EKF, InEKF and IPG observer. The proposed IPG algorithm can obtain the same solution as the default Matlab optimization solver, but in a reduced computational time and a fewer number of iterations. 

%Hence, it can accelerate the optimization step to relieve the computational burden of the MHE approach.

%an observer-based localization approach for nonlinear systems. The proposed method follows the IPG method and can obtain better estimations of the states of the system than the observation values. Two simulation examples and a real-world show the validity of the method.

%% file: sections/appendix.tex
\section*{APPENDIX}
\subsection{Proof of Proposition 1}
This proof mostly follows the proof of Theorem 1 in~\cite{chakrabarti2021accelerating} without the assumption of $\delta = 1$. First, we define the estimation error for the $k$-th iteration as
\begin{align}
    z_k & = \xi_k - \xi^*.
\end{align}
Hence, we can find $z_{k+1}$ using Eq. (\ref{eq:update_x}) as
\begin{align}
    z_{k+1} & = \xi_{k+1} - \xi^* = \xi_k - \delta K_{k} g(\xi_k) - \xi^* \nonumber \\
    & = z_k - \delta K_{k} g(\xi_k).
\end{align}
Using the notation in Eq. (\ref{eq:error_K_k}), we have $K_{k} = \Tilde{K}_{k} + K^*$. From the first order necessary optimality condition, $g(\xi^*) = 0$, hence,
\begin{align}
    & z_{k+1} = z_k - \delta K^* g(\xi_k) - \delta \Tilde{K}_{k} g(\xi_k) \nonumber \\
    & =  \delta K^* (-g(\xi_k) + \frac{1}{\delta} (K^*)^{-1} z_k) - \delta \Tilde{K}_{k} g(\xi_k) \nonumber \\
    & = - \delta K^* \big( g(\xi_k) - g(\xi^*) - \frac{1}{\delta} (K^*)^{-1} z_k \big) -\delta \Tilde{K}_{k} g(\xi_k) \nonumber \\
    & = - \delta K^* \big( g(\xi_k) - g(\xi^*) - \frac{1}{\delta} (H(\xi^*) + \beta I) z_k \big)  -\delta \Tilde{K}_{k} g(\xi_k) \nonumber \\
    & = - \delta K^* \Big( g(\xi_k) - g(\xi^*) - H(\xi^*) z_k - (\frac{1}{\delta}-1) H(\xi^*) z_k \nonumber \\ 
    & \quad - \frac{1}{\delta} \beta z_k \Big) -\delta \Tilde{K}_{k} g(\xi_k) \nonumber \\
    & = - \delta K^* \Big( g(\xi_k) - g(\xi^*) - H(\xi^*) z_k\Big) + \beta K^* z_k \nonumber \\ 
    & \quad +  (1-\delta) K^* H(\xi^*) z_k  -\delta \Tilde{K}_{k} g(\xi_k).
    \label{eq:z_k_plus_1}
\end{align}
We obtain an upper bound on $\| z_{k+1}\|$ as follows.
For the first term in (\ref{eq:z_k_plus_1}), using the fundamental theorem of calculus~\cite{kelley1999iterative},
\begin{align*}
    g(\xi_k) - g(\xi^*) = \int_{0}^{1} H(s \xi_k + (1-s)\xi^*) ds (\xi_k - \xi^*).
\end{align*}
From above we have
\begin{align}
    & g(\xi_k) - g(\xi^*) - H(\xi^*) z_k \nonumber \\
    = &\int_{0}^{1} H(s \xi_k + (1-s)\xi^*) ds  (\xi_k - \xi^*)  - H(\xi^*) z_k \nonumber \\
    % = & (\xi_k - \xi^*) \int_{0}^{1} [H(s \xi_k + (1-s)\xi^*) - H(\xi^*)] ds \nonumber \\
    = & \Big( \int_{0}^{1} [H(s \xi_k + (1-s)\xi^*) - H(\xi^*)] ds \Big) z_k. \label{eqn:grad_hess}
\end{align}
Under Assumption~2, 
\begin{align}
    & \|[H(s \xi_k + (1-s)\xi^*) - H(\xi^{*})]\| \nonumber \\
   & \leq  \gamma \| (s \xi_k + (1-s)\xi^*) - \xi^{*})]\| \nonumber \\
   & = \gamma \| s (\xi_k - \xi^{*}) \| = \gamma  s \|z_k \|. \label{eqn:hess_lip}
\end{align}
By the definition of induced norm,~\eqref{eqn:grad_hess} implies
\begin{align}
    & \| g(\xi_k) - g(\xi^*) - H(\xi^*) z_k \| \nonumber \\
    & \leq \| z_k \| \Big(\int_{0}^{1} \|[H(s \xi_k + (1-s)\xi^*) - H(\xi^{*})]\| ds \Big) \nonumber \\
    & \leq \| z_k \| \Big(\int_{0}^{1} \gamma  s \|z_k \| ds \Big) = \frac{\gamma}{2} \| z_k \|^2, \label{eqn:bd_1}
\end{align}
where the second inequality follows from~\eqref{eqn:hess_lip}.
For the last term in Eq. (\ref{eq:z_k_plus_1}), under Assumption~3 we have
\begin{align}
    & \| \delta \Tilde{K}_{k} g(\xi_k) \| = \| \delta \Tilde{K}_{k} \big( g(\xi_k) - g(\xi^*) \big) \| \nonumber \\
    & \leq  \delta \| \Tilde{K}_{k} \| \|  g(\xi_k) - g(\xi^*)  \| \leq  \delta l \| \Tilde{K}_{k} \| \| \xi_k - \xi^* \| \nonumber \\
    & = \delta l \| \Tilde{K}_{k} \| \| z_k \|. \label{eqn:bd_2}
\end{align}
So, with $\eta = \|K^{*}\|$ and $\| H(\xi^*) \| \leq q$, upon substituting from~\eqref{eqn:bd_1}-\eqref{eqn:bd_2} in~\eqref{eq:z_k_plus_1},
% the upper bound of the estimation error $\|z_{k+1}\|$ is given as,
%\begin{align}
%   & z_{k+1} \nonumber \\
%    & = - \delta K^* \Big( g(\xi_k) - g(\xi^*) - H(\xi^*) z_k\Big) +  (1-\delta) K^* H(\xi^*) z_k \nonumber \\ 
%    & \quad - \beta K^* z_k -\delta \Tilde{K}_{k} g(\xi_k)
%    \label{eq:z_k_plus_2}
%\end{align}
\begin{align}
    & \| z_{k+1}\| \leq  \delta \| K^*\|  \| g(\xi_k) - g(\xi^*) - H(\xi^*) z_k \| \nonumber \\
    & + | 1-\delta | \| K^*\|  \| H(\xi^*) \| \|z_k\| + \beta \| K^*\| \|z_k\| \nonumber \\
    & + \delta \| \Tilde{K}_{k}\|  \|g(\xi_k)\| \nonumber \\
    \leq & \, \frac{\delta \eta \gamma}{2} \| z_k \|^2 + \eta q | 1-\delta | \| z_k \| + \eta \beta \| z_k \| + \delta l \| \Tilde{K}_{k} \| \| z_k \|. \nonumber
\end{align}
Upon substituting $\| \Tilde{K}_{k}\|$ above from~(\ref{eq:K_upper_bound_2}) in Lemma~1,
%If we select $\alpha = \max {\alpha_i}, i=0, \dots, k$, then
%\begin{align}
%    \| z_{k+1}\| \leq & \frac{\eta \gamma}{2} \| z_k \|^2 + (1-\delta) \eta q \|z_k\| + \beta %\eta \|z_k\| \nonumber \\
%    & + l \big( \rho^{k+1} \|\Tilde{K}_0\| + \gamma \eta (\alpha_k \|z_k\|  \nonumber \\
%    & + \rho \alpha_{k-1} \|z_{k-1}\| + \cdots + \rho^k \alpha_{0} \|z_{0}\|) \big)  \| z_k\| 
%\end{align}
\begin{align}
    & \| z_{k+1}\| \leq \frac{\delta \eta \gamma}{2} \| z_k \|^2 + \eta q | 1-\delta | \| z_k \| + \eta \beta \| z_k \| \nonumber \\
    & + \delta l \Big( \rho^{k} \|\Tilde{K}_0\| + \gamma \eta \alpha (\|z_{k-1}\| + \rho \|z_{k-2}\| + \cdots \nonumber \\
    & + \rho^{k-1} \|z_{0}\|) \Big)  \| z_k\|. \label{eqn:zk_bd_1}
\end{align}

Finally, we would prove that $\| z_{k+1}\| < \frac{1}{\mu}\| z_{k}\|$ and $\| z_{k}\| < \frac{1}{\mu \delta \eta \gamma}$ are true for all $k$ using the principle of induction. For $k = 0$, 
%\begin{align}
%    \| z_{k+1}\| \leq &  \frac{1}{\mu}\| z_{k}\| \\
%    \| z_{k}\| < & \frac{1}{\mu \eta \gamma}
%\end{align}
\begin{align}
    \| z_{1}\| \leq &  \| z_0\| (\frac{\delta \eta \gamma}{2} \| z_0 \| + \eta \beta + \eta q | 1-\delta |  + \delta l \| \Tilde{K}_{0}\|).
\end{align}
Hence, if the condition in Eq. (\ref{eq:init_select_cond}) is satisfied,
%\begin{align}
%    \frac{\eta \gamma}{2} \| z_0 \| + \eta \beta + \eta q | 1-\delta | + \delta l \| %\Tilde{K}_{1}\| \leq \frac{1}{2\mu}
%\end{align}
then
\begin{equation*}
    \| z_{1}\| \leq \frac{1}{2\mu}\| z_0\| < \frac{1}{\mu} \| z_0\|.
\end{equation*}
Also, Eq. (\ref{eq:init_select_cond}) implies that $\| z_0\| < \frac{1}{\mu \delta \eta \gamma}$. Therefore, the claims are true for the first iteration.

Next, we suppose that the claims are true for the iteration $1$ to iteration $k$. Then,
\begin{align*}
    \| z_{k+1}\| < \frac{1}{\mu}\| z_{k}\|  < \dots < \frac{1}{\mu^{k+1}} \| z_{0}\| < \frac{1}{\mu^{k+1}} \frac{1}{\mu \delta \eta \gamma}.
\end{align*}
Since $\mu > 1$, the above implies $\| z_{k+1}\| <  \frac{1}{\mu \delta \eta \gamma}$.
In addition, the following inequality is obtained, 
\begin{align*}
    & \|z_k\| + \rho \|z_{k-1}\| + \cdots + \rho^k \|z_{0}\| \nonumber \\
    < & \|z_{0}\| \big( \frac{1}{\mu^{k}} + \frac{\rho}{\mu^{k-1}} + \cdots + \rho^k \big) = \|z_{0}\| \frac{1-(\mu \rho)^{k+1}}{\mu^{k}(1-\mu \rho)}.
\end{align*}
For the iteration $k+1$, in order to show that $\| z_{k+2}\| < \frac{1}{\mu}\| z_{k+1}\|$, from Eq.~(\ref{eqn:zk_bd_1}) and above we have
\begin{align}
    & \| z_{k+2}\| \leq % & \frac{\eta \gamma}{2} \| z_{k+1} \|^2 + \eta \beta \| z_{k+1} \| \nonumber \\
    \frac{\delta \eta \gamma}{2} \| z_{k+1} \|^2 + \eta q | 1-\delta | \| z_{k+1} \| + \eta \beta \| z_{k+1} \| \nonumber \\
    & + \delta l \Big( \rho^{k+1} \|\Tilde{K}_0\| + \gamma \eta \alpha (\| z_{k} \| + \rho \|z_{k-1}\|  \nonumber \\
    & + \rho^2 \|z_{k-2}\| + \cdots + \rho^{k} \|z_{0}\|) \Big)  \| z_{k+1}\| \nonumber \\
    & \leq \| z_{k+1}\| \bigg( \frac{\delta \eta \gamma}{2} \| z_{k+1} \| + \eta q | 1-\delta | + \eta \beta + \delta l \rho^{k+1} \|\Tilde{K}_0\| \nonumber \\
    &  + \delta l \gamma \eta \alpha  \|z_{0}\| \frac{1-(\mu \rho)^{k+1}}{\mu^{k}(1-\mu \rho)} \bigg). \label{eqn:zk2}
\end{align}
%Given $\mu \rho < 1$, and i
If $\alpha$ is selected as 
\begin{equation}
    %\alpha < \frac{\mu^{k}(1-\mu \rho)}{2 \delta l (1-(\mu \rho)^{k+1})}
    \alpha < \frac{\mu^{k}(1-\mu \rho)}{2 l (1-(\mu \rho)^{k+1})},
\end{equation}
then
\begin{equation}
    %\delta l \gamma \eta \alpha  \|z_{0}\| \frac{1-(\mu \rho)^{k+1}}{\mu^{k}(1-\mu \rho)} < \frac{\eta \gamma}{2} \| z_{0} \|
    \delta l \gamma \eta \alpha  \|z_{0}\| \frac{1-(\mu \rho)^{k+1}}{\mu^{k}(1-\mu \rho)} < \frac{\delta \eta \gamma}{2} \| z_{0} \|.
\end{equation}
Since $\rho < 1$, $\delta l \rho^{k+1} \|\Tilde{K}_0\| < \delta l \|\Tilde{K}_0\|$. Then, 
\begin{align}
    & \delta l \gamma \eta \alpha  \|z_{0}\| \frac{1-(\mu \rho)^{k+1}}{\mu^{k}(1-\mu \rho)} + \eta \beta + \eta q | 1-\delta | + \delta l \rho^{k+1} \|\Tilde{K}_0\| \nonumber \\
    & < \frac{\delta \gamma \eta}{2} \| z_{0} \| + \eta \beta + \eta q | 1-\delta | + \delta l \|\Tilde{K}_0\| \leq \frac{1}{2\mu}
\end{align}
Since $\| z_{k+1} \| < \frac{1}{\mu \delta \eta \gamma}$, upon substituting from above in~\eqref{eqn:zk2},
\begin{align}
    \| z_{k+2}\| 
    & < \| z_{k+1}\| (\frac{1}{2\mu} + \frac{1}{2\mu}) = \frac{1}{\mu} \| z_{k+1}\|
\end{align}
Hence, by the principle of induction, we have proved that $\| z_{k+1}\| < \frac{1}{\mu}\| z_{k}\|$ is true for all $k$. As $\mu > 1$, the sequence $\{\| z_k \| = \| \xi_k - \xi^* \|, \forall k \}$ is convergent. 

% \kushal{proof checked up to here}

\subsection{Proof of Proposition 2}
$H(\xi)$ can be expressed as the sum of $N$ matrices, as illustrated in the following form,
\begin{align}
H(\xi) & = {\scriptsize \begin{bmatrix}
* & * &  &  &  & \\
* & * & * &  &  &  \\
    & * & * & * &  &  \\
    &   & \ddots & \ddots & \ddots &  \\
    &  &  &  * & * & * \\
    &  &  &  & * & *
\end{bmatrix}  = \begin{bmatrix}
* & * &  &  &  & \\
* & * & &  &  &  \\
& & &  &  &  \\
&  &  &  &  &   \\
&  &  &  &  &  \\
&  &  &  &  & 
\end{bmatrix}  }\nonumber \\
& {\scriptsize \quad +
\begin{bmatrix}
&  &  &  &  & \\
& * & * &  &  &  \\
& * & * &  &  &  \\
&  &  &  &  &   \\
&  &  &  &  &  \\
&  &  &  &  & 
\end{bmatrix} + \cdots 
 + \begin{bmatrix}
&  &  &  &  & \\
&  &  &  &  &  \\
&  &  &  &  &  \\
&  &  &  &  &   \\
&  &  &  & * & * \\
&  &  &  & * & *
\end{bmatrix} }\nonumber \\
& =: \overline{H}_{0} +  \overline{H}_{1} + \cdots + \overline{H}_{N-1}.
\end{align}
The non-zero part of $\overline{H}_{i}$ is $\overline{\overline{H}}_{i}$ as defined in Eq.~(\ref{eq:bar_bar_H_i}). Moreover, we use the notation $\hat{H}_{(i,j)} \in \mathbb{R}^{n \times n}$ to represent the sub-blocks of $\overline{\overline{H}}_{i}$,
\begin{equation}
\overline{\overline{H}}_{i} = \begin{bmatrix}
\hat{H}_{(i, i)} & \hat{H}_{(i, i+1)} \\
\hat{H}_{(i+1, i)} & \hat{H}_{(i+1, i+1)} 
\end{bmatrix}, \quad (i = 0, \dots, N-1). \nonumber
\end{equation}
Via tedious calculation, we can obtain that 
\begin{align}
    \hat{H}_{(0,0)} = & \, 2\Pi^{-1} + 2 \Tilde{J}_{f(0)} + 2 V_{f(0)}^T \Tilde{Q} \Tilde{H}_{f(0)} + 2 \Tilde{J}_{h(0)} \nonumber \\
    &  + 2 V_{h(0)}^T \Tilde{Q} \Tilde{H}_{h(0)} \nonumber \\
    \hat{H}_{(i, i)} = & \, 2Q^{-1} + 2 \Tilde{J}_{f(i)} + 2 V_{f(i)}^T \Tilde{Q} \Tilde{H}_{f(i)} + 2 \Tilde{J}_{h(i)} \nonumber \\
    &  + 2 V_{h(i)}^T \Tilde{Q} \Tilde{H}_{h(i)}, \quad (i = 1, \dots, N-1) \nonumber \\
    \hat{H}_{(N, N)} = & \, 2Q^{-1} \nonumber \\
    \hat{H}_{(i, i+1)} = & \, -2 J_{f}^T |_{x_i} Q^{-1}, \quad (i = 0, \dots, N-1)\nonumber \\
    \hat{H}_{(i, i-1)} = & \, -2 Q^{-1}J_{f} |_{x_{i-1}}, \quad (i = 1, \dots, N), \nonumber 
\end{align}
with the notations expressed in Section~\ref{sec:convexity}. 
So, if $\overline{\overline{H}}_{i}$ is positive semi-definite, then $\overline{H}_{i}$ is positive semi-definite for all $i = 0, \dots, N-1$. $H(\xi)$ is the sum of $N$ positive semi-definite matrices, which is also positive semi-definite. It concludes that $F(\xi, U, Y)$ at each time instant is convex.